\documentclass[onecolumn]{elsart3p}

\usepackage{amssymb}
\usepackage{wasysym}
\usepackage[english,francais]{babel}

\newtheorem{theorem}{Theorem}[section]

\newtheorem{e-proposition}[theorem]{Proposition}

\newtheorem{e-definition}[theorem]{Definition\rm}
\newtheorem{remark}{\it Remark\/}

\newtheorem{theoreme}{Th\'eor\`eme}[section]

\newtheorem{proposition}[theoreme]{Proposition}

\def\og{\leavevmode\raise.3ex\hbox{$\scriptscriptstyle\langle\!\langle$~}}
\def\fg{\leavevmode\raise.3ex\hbox{~$\!\scriptscriptstyle\,\rangle\!\rangle$}}
\journal{the Acad\'emie des sciences}
\begin{document}
\centerline{Mathematical Problems in Mechanics}
\begin{frontmatter}


\title{$3D-2D$ dimensional reduction for a nonlinear optimal design problem with perimeter penalization}


\selectlanguage{english}
\author[label1]{Gra\c{c}a Carita},
\ead{gcarita@uevora.pt}
\author[label2]{Elvira Zappale}
\ead{ezappale@unisa.it}

\address[label1]{CIMA-UE, Departamento de Matem\'{a}tica,
Universidade de \'Evora, Rua Rom\~{a}o Ramalho, 59 7000-671 \'Evora, Portugal.}
\address[label2]{D.I.E.I.I., Universit\`a degli Studi di Salerno, Via Ponte Don Melillo, 84084 Fisciano (SA)
Italy.}



\begin{abstract}
\noindent A 3D-2D dimension reduction for a nonlinear optimal design problem with a
perimeter penalization is performed in the realm of $\Gamma$-convergence,
providing an integral representation for the limit functional.

\vskip 0.5\baselineskip

\selectlanguage{francais}
\noindent{\bf
R\'esum\'e} \vskip 0.5\baselineskip \noindent{\bf R\'eduction dimensionnelle $3D-2D$
d'un probl\`eme non lin\'eaire d'optimisation de forme avec p\'enalisation
sur le p\'erim\`etre} On effectue dans ce travail une r\'eduction dimensionnelle
3D-2D d'un probl\`{e}me non lin\'{e}aire d'optimisation de forme avec une
p\'{e}nalisation du p\'{e}rim\`{e}tre. Une repr\'esentation int\'egrale de la fonctionnelle limite est obtenue.
\end{abstract}
\end{frontmatter}

\selectlanguage{francais}
\section*{Version fran\c{c}aise abr\'eg\'ee}

\noindent On s'int\'eresse dans ce travail au comportement asymptotique
d'une suite de probl\`{e}mes non lin\'{e}aires d'optimisation de forme avec p\'{e}nalisation du p\'{e}rim\`{e}tre sur le domaine cylindrique 
$\Omega \left( \varepsilon \right) :=$ $\omega \times \left( -\varepsilon
,\varepsilon \right) ,$ 
o\`{u} $\varepsilon >0$ et $\omega $ est un ouvert
born\'{e} de $\mathbb{R}^{2}.$ On suppose que le domaine $\Omega \left( \varepsilon
\right) $, occup\'e par le solide, est constitu\'e par deux mat\'{e}riaux hyper\'{e}lastiques dont les
densit\'{e}s d'\'{e}nergie, $W_{i}:\mathbb{R}^{3\times 3}\rightarrow \mathbb{%
R}$, $i=1,2,$ sont continues, v\'erifiant la condition de
croissance 
\begin{equation}
\beta ^{\prime }\left( \left\vert F \right\vert ^{p}-1\right) \leq
W_{i}(F)\leq \beta (1+|F |^{p})\hbox{ pour tout }F \in \mathbb{R}%
^{3\times 3},\;\;p>1,\;\ i=1,2,\hbox{ avec }\beta \geq \beta
^{\prime }>0. \label{H1}
\end{equation}
Plus pr\'ecis\'ement, on consid\`{e}re le probl\`{e}me de minimisation suivant
\begin{equation}
\begin{array}{l}
\inf_{%
\begin{array}{c}
v\in W^{1,p}(\Omega (\varepsilon );\mathbb{R}^{3}) \\
\chi _{E(\varepsilon )}\in BV(\Omega (\varepsilon );\{0,1\})%
\end{array}%
}\Big\{\frac{1}{\varepsilon }\Big(\int_{\Omega (\varepsilon )}\left( \chi
_{E(\varepsilon )}W_{1}+(1-\chi _{E(\varepsilon )})W_{2}\right) (\nabla
v)dx-\int_{\Omega (\varepsilon )}{{f}_\varepsilon}\cdot v~dx \\
\;\;\;\;\;\;\;\;\;\;\;\;\;\;\;\;\;\,\,\,\;\;\;\;\;\;\;\;\;\;\;\;\;\;\;\;\;\;%
\;\;\;\;\;+P(E(\varepsilon );\Omega (\varepsilon ))\Big):\;\;v=0\hbox{ sur }%
\partial \omega \times (-\varepsilon ,\varepsilon ),\frac{1}{\mathcal{L}%
^{3}(\Omega (\varepsilon ))}\int_{\Omega (\varepsilon )}\chi _{E(\varepsilon
)}~dx=\lambda \Big\},%
\end{array}
\label{originalpb}
\end{equation}%
o\`{u} $E\left( \varepsilon \right) \subset \Omega \left( \varepsilon
\right) $ est un ensemble mesurable de p\'{e}rim\`{e}tre fini et ${%
f_{\varepsilon}}\in L^{p^{\prime }}\left( \Omega (\varepsilon );\mathbb{R}%
^{3}\right) $ o\`{u} $\frac{1}{p}+\frac{1}{p^{\prime }}=1$ et $\lambda \in [0,1]$ est la fraction du volume rempli par le premier mat\'{e}riau.

Nous commen\c{c}ons par effectuer un changement de variables afin de rendre le domaine ind\'ependant de $\varepsilon$ (cf. (\ref{1edilation})). On obtient ainsi le probl\`{e}me suivant%
\begin{equation}
\begin{array}{r}
\inf\limits_{%
\begin{array}{l}
u\in W^{1,p}(\Omega ;\mathbb{R}^{3}) \\
\chi \in BV(\Omega ;\{0,1\})%
\end{array}%
}\Big\{\int_{\Omega }(\chi W_{1}+(1-\chi )W_{2})\left( \nabla _{\alpha }u|%
\frac{1}{\varepsilon }\nabla _{3}u\right) dx-\int_{\Omega }f\cdot
udx+\left\vert \Big(D_{\alpha }\chi |\frac{1}{\varepsilon }D_{3}\chi \Big)%
\right\vert (\Omega ): \\
\qquad \qquad \qquad \qquad \qquad \qquad \qquad \left. u=0\hbox{ sur }%
\partial \omega \times (-1,1),\frac{1}{\mathcal{L}^{3}(\Omega )}\int_{\Omega
}\chi dx=\lambda \Big\}\right. .%
\end{array}
\label{P1}
\end{equation}%
En utilisant les techniques de $\Gamma -$convergence, on d\'emontre que, lorsque $\varepsilon $
tend vers z\'ero, le probl\`{e}me (\ref{P1}) converge vers un probl\`{e}me nonlin\'{e}aire bidimensionnel.
 
\noindent Il suit un r\'esultat g\'{e}n\'{e}ral o\`{u} la p\'{e}nalisation du p\'{e}rim\`{e}tre dans le probl\`eme initial est remplac\'{e}e par une int\'{e}grale
elliptique. Ainsi, on consid\`{e}re $\Psi :\mathbb{R}%
^{3}\rightarrow \lbrack 0,+\infty \lbrack $ une fonction paire, continue,
positivement homog\`{e}ne de degr\'{e} $1$ et telle que
\begin{equation}
\displaystyle{\exists C\in ]0,+\infty \lbrack \,\,:\,\,\forall \nu \in
\mathbb{R}^{3}\quad \frac{1}{C}|\nu |\leq \Psi (\nu )\leq C|\nu |.}
\label{Psigrowth}
\end{equation}%

On \'{e}tudie donc,  le comportement asymptotique, lorsque $\varepsilon $ tend vers
z\'ero, du probl\`{e}me suivant%
\begin{equation}
\begin{array}{ll}
\inf_{%
\begin{array}{c}
v\in W^{1,p}(\Omega (\varepsilon );\mathbb{R}^{3}) \\
\chi \in BV(\Omega (\varepsilon );\{0,1\})%
\end{array}%
}\Big\{\frac{1}{\varepsilon }\Big(\int_{\Omega (\varepsilon )}\Big(\chi
_{E(\varepsilon )}W_{1}+(1-\chi _{E(\varepsilon )})W_{2}\Big)(\nabla
v)dx-\int_{\Omega (\varepsilon )}{f}_{\varepsilon}\cdot v~dx & \\
\quad \quad \quad \quad \quad \quad \quad \quad \quad \quad +\int_{\partial
E(\varepsilon )}\Psi (\nu _{E(\varepsilon )})d\mathcal{H}^{2}\Big ):v=0%
\hbox{ sur }\partial \omega \times (-\varepsilon ,\varepsilon ),\frac{1}{%
\mathcal{L}^{3}(\Omega (\varepsilon ))}\int_{\Omega (\varepsilon )}\chi
_{E(\varepsilon )}~dx=\lambda \Big\}, &
\end{array}
\label{originalpb2}
\end{equation}%
o\`{u} $\mathcal{H}^{2}$ d\'{e}signe la mesure de Hausdorff restreinte \`a $\partial E(\varepsilon )$ et $\nu _{E(\varepsilon )}$ est la normale ext%
\'{e}rieure \`{a} $ E(\varepsilon )$.
On obtient finalement, le probl\`{e}me limite (\ref{limitpb2}).
\selectlanguage{english}
\section{Introduction and setting of the problem}
\label{introduction}

The study of thin structures has been the object of many investigations. In particular, in mechanical engineering it is important for applications to minimize, under a given system of loads, the compliance (namely, the
opposite of the total energy at equilibrium) of a given structure, satisfying
a constraint on the volume. In order to design thin structures with the best possible resistance-weight ratio, the asymptotic behaviour of the compliance as the thickness of the sample tends to zero is studied. For a background on the modelling of thin plates we refer to the monographs of \cite{A} and \cite{C}. 

 Let $\Omega\left(  \varepsilon\right)  :=\omega\times\left(
-\varepsilon,\varepsilon\right),$ where $\omega$ is a bounded open domain of
$\mathbb{R}^{2}$ and $\varepsilon>0$, and for the sake of illustration let us assume that $\Omega
(\varepsilon)$ is clamped on its lateral boundary. We suppose also
that $\Omega(\varepsilon)$ is filled with two  materials with
respective energy densities $W_{1}$ and $W_{2}$, where $W_{i}:\mathbb{R}%
^{3\times3}\rightarrow\mathbb{R}$, $i=1,2$ are continuous functions satisfying
(\ref{H1}) (see Remark 2 below, where assumption (\ref{H1}) is discussed).


Let $E(\varepsilon)$ be the first phase, $f_{\varepsilon}$ be the given load on
$\Omega(\varepsilon)$ and assume that the volume fraction of each phase is
given by $\lambda:=\frac{1}{\mathcal{L}^{3}(\Omega(\varepsilon))}\int%
_{\Omega(\varepsilon)}\chi_{E(\varepsilon)}(x)dx\in\lbrack0,1]$, where
$\chi_{E(\varepsilon)}$ denotes the characteristic function of the phase
$E(\varepsilon)$. The compliance
$C^{\varepsilon}(\chi_{E(\varepsilon)})$ is defined as
$$
C^{\varepsilon}(\chi_{E(\varepsilon)})  :=-\inf_{v\in
W^{1,p}(\Omega(\varepsilon);\mathbb{R}^{3})}\Big\{  {\frac{1}{\varepsilon
}\int_{\Omega(\varepsilon)}\left(  \left(  \chi_{E(\varepsilon)}W_{1}%
+(1-\chi_{E(\varepsilon)})W_{2}\right)  (\nabla v)-{f_{\varepsilon}}\cdot
v\right)  dx:}
v=0\hbox{ on }\partial\omega\times(-\varepsilon,\varepsilon)\Big\}.
$$

In  \cite{BFS1}, \cite{BFS2} and \cite{BFS3} the asymptotic behaviour of a $3D$ optimal elastic compliance problem is studied, as the thickness (or the cross section in the case of beams) tends to zero and the volume fraction in the design region remains unchanged. It is assumed that the material has a convex and $2$-homogeneous potential and the analysis is performed in the small-displacement setting.
The prescription of the volume in the minimum problem can be dropped  by adding a Lagrange multiplier to penalize the volume in the cost functional.  The asymptotic analysis performed in these papers leads to a fictitious material with local density, taking all the values in $[0,1]$ and not to a precise limit set, due to the loss of compactness in the characteristic functions.

In this work we focus our attention on studying  the worst possible design of a two-phase mixture of elastic materials in a thin film in the same spirit of \cite{FF} and \cite{BFF}, where the asymptotic analysis of a two-field minimization problem has been studied (i.e. $(\chi,
u)$, (design region, deformation)) as the thickness of the sample tends to zero. Having in mind
the results contained in \cite{AB} and \cite{LK}, we introduce a perimeter
penalization in our functional in order to derive from the $3D$ energy a
limiting $2D$ model where the design region is explicitly determined, and we refer to \cite{CFP} for a detailed study about regularity of the limits (set, deformation).




Let us consider the optimal design problem in (\ref{originalpb}) where $E\left(
\varepsilon\right)  \subset\Omega\left(  \varepsilon\right)  $ is a measurable
subset of $\Omega(\varepsilon)$ with finite perimeter, i.e.,
\begin{equation}
P(  E(  \varepsilon)  ;\Omega(  \varepsilon))  :=\sup\Big\{  \int_{E(\varepsilon)}{\it div}%
\varphi\,dx:\varphi\in C_{c}^{1}(\Omega(\varepsilon);\mathbb{R}^{3}%
),\Vert\varphi\Vert_{L^{\infty}}\leq1\Big\}  <+\infty, \label{perimeter1}%
\end{equation}
and we assume that the load ${f_{\varepsilon}}\in L^{p^{\prime}}(\Omega(\varepsilon
);\mathbb{R}^{3})$, with $\frac{1}{p}+\frac{1}{p^{\prime}}=1$.

In order to study the asymptotic behavior of (\ref{originalpb}) we reformulate our problem in a fixed 3D domain through a $\frac{1}{\varepsilon}$- dilation in the transverse direction $x_3$ and then we perform
$\Gamma-$ convergence with respect to the pair (design region, deformation).  Set
\begin{equation}%
\begin{array}[c]{ll}%
\Omega:=\omega\times\left(  -1,1\right)  ,\quad
E_{\varepsilon}:=\left\{  \left(  x_{1},x_{2},x_{3}\right)  \in\Omega:\left(
x_{1},x_{2},\varepsilon x_{3}\right)  \in E\left(  \varepsilon\right)
\right\}  ,\\
u\left(  x_{1},x_{2},x_{3}\right)  :=v\left(  x_{1},x_{2},\varepsilon
x_{3}\right)  ,\quad
f\left(  x_{1},x_{2},x_{3}\right)  :={f_{\varepsilon}}\left(  x_{1}%
,x_{2},\varepsilon x_{3}\right)  , \quad
\chi_{E_{\varepsilon}}\left(  x_{1},x_{2},x_{3}\right)  :=\chi_{E\left(
\varepsilon\right)  }\left(  x_{1},x_{2},\varepsilon x_{3}\right)  ,
\end{array}
\label{1edilation}%
\end{equation}
where $v$ is any admissible field for (\ref{originalpb}).

 In the sequel we will denote
$x_{\alpha}:=(x_{1},x_{2}),$ $dx_{\alpha}:=dx_{1}dx_{2}$ and $\nabla_{\alpha}$
and $D_{\alpha}$ will be identified with the pair $\left(  \nabla_{1}%
,\nabla_{2}\right)  ,~\left(  D_{1},D_{2}\right)  ,$ respectively. For every
matrix $\overline{F}\in\mathbb{R}^{3\times2}$ and any $z\in\mathbb{R}^{3}$,
$F:=(\overline{F}|z)$ represents the matrix in $\mathbb{R}^{3\times3}$ whose
first two columns are those of $\overline{F}$ and the last column is given by
the vector $z$.

Observe that by (\ref{perimeter1}) and using the definition of total variation,
$P\left(  E\left(  \varepsilon\right)  ;\Omega\left(  \varepsilon\right)
\right)  =\left\vert D\chi_{E\left(  \varepsilon\right)  }\right\vert \left(
\Omega\left(  \varepsilon\right)  \right)  $.
Making the change
of variables $y_{3}:=\varepsilon x_{3}$ and $y_{\alpha}:=x_{\alpha}$ we have $\frac{1}{\varepsilon}\left\vert D\chi_{E(\varepsilon)}\right\vert \left(
\Omega(\varepsilon)\right)  =\left\vert \left(  D_{\alpha}\chi_{\varepsilon
}\left|\frac{1}{\varepsilon}D_{3}\chi_{\varepsilon}\right.\right)  \right\vert (\Omega),$ where $\chi_\varepsilon:=\chi_{E_{\varepsilon}}$ stands for the
characteristic function of $E_{\varepsilon}$. Hence
%
we are led to the rescaled minimum problem (\ref{P1}).



For every $\varepsilon>0$, let $J_{\varepsilon}:L^{1}(\Omega;\{0,1\})\times
L^{p}(\Omega;\mathbb{R}^{3})\rightarrow\lbrack0,+\infty]$ be the functional
defined as follows
\begin{equation}
J_{\varepsilon}(\chi,u):=\left\{
\begin{array}
[c]{l}%
{%
\int_{\Omega}
\left(  \chi W_{1}\left(  \nabla_{\alpha}u \left| \frac{1}{\varepsilon}\nabla
_{3}u\right. \right)  +(1-\chi)W_{2}\left(  \nabla_{\alpha}u\left|\frac{1}{\varepsilon
}\nabla_{3}u\right.\right)  \right)  dx}-{%
\displaystyle\int_{\Omega}
f\cdot udx+\left\vert\left( D_{\alpha}\chi\left|\frac{1}{\varepsilon}D_{3}\chi\right)\right.\right\vert
(\Omega)}\\
\;\;\;\,\,\,\;\;\;\;\;\;\;\;\;\;\;\;\,\,\,\;\;\;\;\;\;\;\;\;\;\;\;\,\,\,\;\;\;\;\;\;\;\;\;\;\;\;\;\;\;\;\;\;\;\ \ \ \ \ \ \ \ \ \ \ \ \ \ \ \ \ \ \ \ \ \ \ \ \ \hbox{ in }BV(\Omega
;\{0,1\})\times W^{1,p}(\Omega;\mathbb{R}^{3}),\\
+\infty
\;\;\;\,\,\,\;\;\;\;\;\;\;\;\;\;\;\;\,\,\,\;\;\;\;\;\;\;\;\;\;\;\;\;\;\;\;\;\,\,\,\;\;\;\;\;\;\;\;\ \ \ \ \ \ \ \ \ \ \ \ \ \ \ \ \ \ \ \ \ \ \ \ \ \;\hbox{otherwise.}
\end{array}
\right.  \label{Je}%
\end{equation}

Let $V:\{0,1\} \times \mathbb R^{3 \times 3} \rightarrow \lbrack 0,+\infty)$ be given by 
\begin{equation}\label{V}
\displaystyle{V(\chi, F):=\chi W_1(F)+ (1-\chi)W_2(F),}
\end{equation}
with $W_1$ and $W_2$ satisfying (\ref{H1}).
Analogously, let $\overline{V}:\{0,1\}\times\mathbb{R}^{3\times2}\rightarrow\lbrack
0,+\infty)$ be defined as
\begin{equation}
\overline{V}\left(  \chi,\overline{F}\right)  :=\chi\overline{W}_{1}\left(
\overline{F}\right)  +\left(  1-\chi\right)  \overline{W}_{2}\left(
\overline{F}\right)  , \label{Vbar}%
\hbox{ with }
\overline{W}_{i}\left(  \overline{F}\right)  :=\inf_{c\in\mathbb{R}^{3}}%
W_{i}\left(  \overline{F}|c\right)  ,\qquad\overline{F}\in\mathbb{R}%
^{3\times2},\;\;\;\;i=1,2.
\end{equation}

Consider the functional $J_{0}:L^{1}(\Omega;\{0,1\})\times L^{p}%
(\Omega;\mathbb{R}^{3})\rightarrow\lbrack0,+\infty]$ as
\begin{equation}
J_{0}(\chi,u):=\left\{
\begin{array}
[c]{ll}%
{2%
{\displaystyle\int_{\omega}}
Q\overline{V}(\chi,\nabla_{\alpha}u)dx_{\alpha}-%
{\displaystyle\int_{-1}^{1}}
{\displaystyle\int_{\omega}}
f\cdot udx_{\alpha}dx_{3}+2|D_\alpha \chi|(\omega),} & \hbox{if }(\chi,u)\in
BV(\omega;\{0,1\})\times W^{1,p}(\omega;\mathbb{R}^{3}),\medskip\\
+\infty & \hbox{otherwise,}%
\end{array}
\right.  \label{J0}%
\end{equation}
where $Q\overline{V}$ stands for the quasiconvexification of $\overline{V}$ in
the second variable. Namely, for every $(\chi,\overline{F})\in\{0,1\}\times\mathbb{R}^{3\times2}$
\begin{equation}
Q\overline{V}\left(  \chi,\overline{F}\right)  :=\inf\left\{  \int_{Q^{\prime
}}\overline{V}\left(  \chi,\overline{F}+\nabla_{\alpha}\varphi\right)  dx_{\alpha
}:\varphi\in W_{0}^{1,p}\left(  Q^{\prime};\mathbb{R}^{3}\right)  \right\}  , 
\label{QVbar}%
\end{equation}
\noindent where  $Q^{\prime} \subset \mathbb R^2$ denotes the unit cube.

\noindent Our main result is the following.

\begin{theorem}
\label{mainthm}  The family of functionals $\{J_{\varepsilon}\}$ $\Gamma
$-converges, with respect to the strong topology of $L^{1}(\Omega
;\{0,1\})\times L^{p}(\Omega;\mathbb{R}^{3})$ to $J_{0}$, as
$\varepsilon\rightarrow0^{+}$.
\end{theorem}

\smallskip

\begin{remark}
We observe that Theorem \ref{mainthm} entails the convergence, as
$\varepsilon\rightarrow0^{+}$, of problems (\ref{originalpb}) in their
rescaled version (\ref{P1}) to the problem
\begin{equation}
\inf\limits_{%
\begin{array}
[c]{l}%
u\in W_{0}^{1,p}(\omega;\mathbb{R}^{3})\\
\chi\in BV(\omega;\{0,1\})
\end{array}
}\left\{  2\int_{\omega}Q\overline{V}(\chi,\nabla_{\alpha}u)dx_{\alpha}%
-\int_{-1}^{1}\int_{\omega}f\cdot udx_\alpha dx_3+2|D_\alpha\chi|(\omega):\frac{1}{\mathcal{L}%
^{2}(\omega)}\int_{\omega}\chi dx_{\alpha}=\frac{1}{2}\lambda\right\}  .
\label{limitPb1}%
\end{equation}

\noindent In fact, this is due to the strong convergence in $L^{1}(\Omega;\{0,1\})\times
L^{p}(\Omega;\mathbb{R}^{3})$, of the sequence of almost minimizers
$\{(\chi_{\varepsilon}, u_{\varepsilon})\}$ of  (\ref{P1}) to $(\chi, u)
\in BV(\omega;\{0,1\})\times W^{1,p}_{0}(\omega;\mathbb{R}^{3})$.
And so, the constraint in the volume fraction $\frac{1}{\mathcal{L}^{3}%
(\Omega(\varepsilon))}\int_{\Omega(\varepsilon)}\chi_{E(\varepsilon)} dx =
\frac{1}{\mathcal{L}^{3}(\Omega)}\int_{\Omega}\chi_{\varepsilon}dx =\lambda$
is kept in the limit, as well as the boundary conditions (cf. Remark
\ref{compactness}).

\end{remark}

It is worthwhile to compare Theorem \ref{mainthm} with similar results in
\cite{FF} and \cite{BFF}. To this aim we observe that, in spite of what is proven
therein, namely $\Gamma$-convergence results with respect to the convergence
$L^{\infty}_{\mathrm{weak}\ast}\times L^{p}$ for $(\chi, u)$, the presence of
the perimeter in our energy (\ref{Je}), allows us to have a stronger
convergence on the characteristic functions and thus to determine the worst possible design set.
On the other hand, the fact that the perimeter is inserted in our model leads naturally to compare our results with those contained in \cite{AB}.  Indeed, if  $W_{1}$ and $W_{2}$ are of type
$W_{1}\left(  \cdot\right)  :=\alpha^{\prime}\left\vert \cdot\right\vert ^{2}$
and $W_{2}\left(  \cdot\right)  :=\alpha\left\vert \cdot\right\vert ^{2}$, with
$0<\alpha^{\prime}<\alpha$ suitable constants, clearly $Q\overline{V}(\chi,\overline{F})$ coincides
with $\alpha^{\prime}\chi|\overline{F}|^{2}+\alpha(1-\chi)|\overline{F}|^{2}$.
Hence, by \cite[Theorem 2.2]{AB} the solution of the minimum problem
(\ref{limitPb1}) is locally H\"{o}lder continuous and the optimal design set
is equivalent to an open set $A\times(-1,1)$, $A\subset\omega$. More refined
results, about regularity in $2D$, in the convex setting, may be found in
\cite{L1,L2} and in the references quoted therein. 

In our model, the lack of convexity in $W_1$ and $W_2$ entails, as underlined by Proposition \ref{continuityQVbar}, that we obtain a
limit energy which depends continuously on the characteristic function of the design set and requires a quasiconvexification procedure in the deformation variable. We refer to \cite{CFP} for regularity results related to our setting.


Details about the results are contained in the next section, while for the
properties related to $\Gamma$-convergence, sets of finite perimeter and $BV$
functions we refer to \cite{DM} and \cite{AFP2}, respectively.

\section{The limit problem}

We start by stating the properties of the energy densities in (\ref{V}) and (\ref{Vbar}) that we will exploit in the sequel.

\begin{proposition}
\label{continuityVbar}
 Let $\overline{V}$ be as in $\left(  \ref{Vbar}\right)  .$ 
Then $\overline{V}$ is continuous and
satisfies%
\begin{equation}
\displaystyle{\beta^{\prime}\left(\left\vert \overline{F}\right\vert ^{p}-1\right)\leq\overline{V}\left(
\chi,\overline{F}\right)  \leq\beta\left(  1+\left\vert \overline
{F}\right\vert ^{p}\right)  ,} \label{growth}%
\end{equation}
where $\beta^{\prime}$ and $\beta$ are the constants in $\left(
\ref{H1}\right)  .$ Moreover, $ \left|  \overline{V}\left(  \chi,\overline{F}\right)
-\overline{V}(\chi^{\prime}, \overline{F})\right|  \leq2 \beta\left|
\chi-\chi^{\prime}\right|  (1+ |\overline{F}|^{p}) .$

\end{proposition}

\begin{proposition}
\label{continuityQVbar} 
The  function $Q\overline{V}$ in (\ref{QVbar})  is  continuous and satisfies (\ref{growth}), and 
\begin{equation}
\label{ucQVbar}
\displaystyle{|Q\overline{V}(\chi, \overline{F})-Q\overline{V}(\chi^{\prime},
\overline{F})|\leq C|\chi^{\prime}- \chi|(1+|\overline{F}|^{p}).}
\end{equation}
\end{proposition}

\begin{remark}
\label{compactness} We claim that energy bounded sequences  $\left\{  \left(  \chi_{\varepsilon},u_{\varepsilon}\right)
\right\}  $ for problem
(\ref{P1}), with $u_{\varepsilon}$ clamped on $\partial\omega\times(-1,1)$,
are compact in $L^{1}(\Omega;\{0,1\})\times L^{p}(\Omega;\mathbb{R}^{3})$ and
with limit in $L^{1}(\omega;\{0,1\})\times L^{p}(\omega;\mathbb{R}^{3})$.

\noindent If $\left\{  (\chi_{\varepsilon},u_{\varepsilon})\right\}  $ is a sequence
such that $J_{\varepsilon}(\chi_{\varepsilon},u_{\varepsilon})\leq C$, then there exists
$C^{\prime} \in \mathbb R^+$ such that the following bounds hold
$$
\Vert u_{\varepsilon}\Vert_{W^{1,p}}\leq C^{\prime},\qquad\left\Vert \frac
{1}{\varepsilon}\nabla_{3}u_{\varepsilon}\right\Vert _{L^{p}}\leq C^{\prime}, \;\;\;
\left\vert \left(  D_{\alpha}\chi_{\varepsilon}\left\vert \frac{1}%
{\varepsilon}D_{3}\chi_{\varepsilon}\right.  \right)  \right\vert (\Omega)\leq
C^{\prime}.
$$
An argument entirely similar to that exploited in \cite[Lemma 3]{LDR}, entails that there exists $u\in W^{1,p}(\Omega;\mathbb{R}^{3})$ such
that $\nabla_{3}u\equiv 0$, and so $u$ can be identified with a function (still
denoted in the same way) $u\in W^{1,p}(\omega;\mathbb{R}^{3}).$ Thus we may find a subsequence, not relabelled, $\{u_\varepsilon\}$ such that $u_{\varepsilon}\rightharpoonup u\hbox{ in }W^{1,p}(\Omega;\mathbb{R}^{3})$, and a measurable set $E\subset\Omega$ such that $\chi_{\varepsilon
}{\rightharpoonup}\ast\chi_{E}$ and $D_{3}\chi_{E}\equiv0.$
Hence, there exists $E^{\prime}\subset\omega$, with $|D\chi_{E}|(\Omega)=2|D\chi_{E^{\prime}}|(\omega), $ where $E=E^{\prime}\times(-1,1).$ In the following we will identify the set $E$  with the set $E^{\prime}$ and denote
$\chi_{E^{\prime}}$ by $\chi$. 

\smallskip

\end{remark}

\smallskip

\noindent {\bf Proof of Theorem \ref{mainthm}} For every $\varepsilon>0$, let $J_\varepsilon$ be the functional in (\ref{Je}). The $\Gamma$-convergence with respect to the
separable metric space $L^{1}(\Omega;\{0,1\})\times L^{p}(\Omega
;\mathbb{R}^{3})$ ensures that for each sequence $\{\varepsilon\}$ there
exists a subsequence, still denoted by $\{\varepsilon\}$, such that
$\Gamma-\lim_{\varepsilon\rightarrow0^{+}}(L^{1}(\Omega;\{0,1\})\times
L^{p}(\Omega;\mathbb{R}^{3}))J_{\varepsilon}$ exists.

\noindent For every $(\chi,u)\in L^{1}(\Omega;\{0,1\})\times L^{p}(\Omega;\mathbb{R}%
^{3})$, let $J\left(  \chi,u\right)$ be its  $\Gamma$-limit. By virtue of Urysohn property, it suffices  to prove that
any sequence $\{J_{\varepsilon}\}$ admits a further subsequence whose $\Gamma
$-limit, $J(\chi,u),$ coincides with $J_{0}(\chi,u)$ in (\ref{J0}).

We observe that if $(\chi,u)\in(L^{1}(\Omega;\{0,1\})\times L^{p}%
(\Omega;\mathbb{R}^{3}))\setminus(BV(\omega;\left\{  0,1\right\}  )\times
W^{1,p}(\omega;\mathbb{R}^{3}))$, then $J(\chi,u)=+\infty$. Indeed, if this is
not the case, from $J(\chi,u)<+\infty$ we would get the existence of a
sequence $\{(\chi_{\varepsilon},u_{\varepsilon})\}$ converging to $(\chi,u)$
such that $J_{\varepsilon}(\chi_{\varepsilon},u_{\varepsilon})<+\infty$ and by
Remark \ref{compactness} this would imply $(\chi,u)\in BV(\omega;\left\{
0,1\right\}  )\times W^{1,p}(\omega;\mathbb{R}^{3}).$

The remaining proof is divided into two steps. First we show the lower bound,
then we prove the upper bound.

\noindent\textsl{Step one:} We claim that for every $(\chi,u)\in
BV(\omega;\{0,1\})\times W^{1,p}(\omega;\mathbb{R}^{3})$
 
$$
J(\chi,u)\geq2\int_{\omega}Q\overline{V}\left(  \chi\left(  x_{\alpha}\right)
,\nabla_{\alpha}u\left(  x_{\alpha}\right)  \right)  dx_{\alpha}-\int_{-1}%
^{1}\int_{\omega}f\left(  x_{\alpha},x_{3}\right)  u\left(  x_{\alpha}\right)
dx_{\alpha}dx_{3}+2\left\vert D_\alpha\chi\right\vert \left(  \omega\right)  .
$$
To prove the claim, let $\{(\chi_{\varepsilon},u_{\varepsilon})\}\subset
L^{1}(\Omega;\{0,1\})\times L^{p}(\Omega;\mathbb{R}^{3})$ be a sequence
converging to $(\chi,u)\in BV(\omega;\{0,1\})\times W^{1,p}(\omega
;\mathbb{R}^{3})$. For the forces and for the perimeter the lower bound
follows by $L^{p}$ strong convergence of $\left\{  u_{\varepsilon}\right\}  $
and lower semicontinuity of the perimeter, respectively.

For what concerns the bulk energy, by virtue of the Decomposition
Lemma for scaled gradients, (cf. \cite[Theorem 1.1]{BF}) there exist a subsequence of $\{u_\varepsilon\}$, not relabelled, and a sequence $\{w_{\varepsilon}\}$ converging to $u\in W^{1,p}(\omega;\mathbb{R}^{3})$, such that the scaled
gradients $\{\left(  \nabla_{\alpha}w_{\varepsilon},\frac{1}{\varepsilon
}\nabla_{3}w_{\varepsilon}\right)  \}$ are $p$-equiintegrable, and ${\cal L}^3(\Omega \setminus A_\varepsilon)\to 0$ as $\varepsilon \to 0^+$, where  $A_\varepsilon :=\{x \in \Omega: u_{\varepsilon}\equiv
w_{\varepsilon}\}$. Denoting the bulk energy density of $J_\varepsilon$ by $V$ as in (\ref{V}), one obtains
\begin{equation}\label{chainofineq}
\begin{array}{ll}
&  \liminf_{\varepsilon\rightarrow0^{+}}\int_{\Omega}V\left(  \chi
_{\varepsilon},\left(  \nabla_{\alpha}u_{\varepsilon}\left|\frac{1}%
{\varepsilon}\nabla_{3}u_{\varepsilon}\right.\right) \right)  dx\\
&  \geq\liminf_{\varepsilon\rightarrow 0^{+}}\int_{\Omega}V\left(
\chi_{\varepsilon},\left(  \nabla_{\alpha}w_{\varepsilon}\left|\frac
{1}{\varepsilon}\nabla_{3}w_{\varepsilon}\right.\right)\right)  dx
 -\beta\limsup_{\varepsilon\rightarrow0^{+}}\int_{\Omega\setminus
A_{\varepsilon}}\left(  1+\left\vert \left(  \nabla_{\alpha}w_{\varepsilon
}\left|\frac{1}{\varepsilon}\nabla_{3}w_{\varepsilon}\right.\right)  \right\vert
^{p}\right)  dx\\
&  \geq\liminf_{\varepsilon\rightarrow0^{+}}\int_{\Omega}V\left(
\chi_{\varepsilon},\left(  \nabla_{\alpha}w_{\varepsilon}\left|\frac
{1}{\varepsilon}\nabla_{3}w_{\varepsilon}\right.\right) \right)  dx
 \geq\liminf_{\varepsilon\rightarrow0^{+}}\int_{\Omega}\overline{V}\left(
\chi_{\varepsilon},\nabla_{\alpha}w_{\varepsilon}\right)  dx\geq
\liminf_{\varepsilon\rightarrow0^{+}}\int_{\Omega}Q\overline{V}\left(
\chi_{\varepsilon},\nabla_{\alpha}w_{\varepsilon}\right)  dx.
\end{array}
\end{equation}

Observe that, by (\ref{ucQVbar})
\begin{equation}\label{marker}
\int_{\Omega}|Q\overline{V}\left(  \chi_{\varepsilon},\nabla_{\alpha
}w_{\varepsilon}\right)  -Q\overline{V}\left(  \chi,\nabla_{\alpha
}w_{\varepsilon}\right)  |\,dx\leq C\int_{\Omega}|\chi_{\varepsilon}%
-\chi|(1+|\nabla_{\alpha}w_{\varepsilon}|^{p})\,dx.
\end{equation}
Thus, the $p$-equiintegrability of $\left\{  \left(  \nabla_{\alpha
}w_{\varepsilon}\left\vert \frac{1}{\varepsilon}\nabla_{3}w_{\varepsilon
}\right.  \right)  \right\}  $ and (\ref{marker}) ensure that as $\varepsilon
\rightarrow0^{+}$, $\chi_\varepsilon$ can be replaced by $\chi$ in the right-hand side of (\ref{chainofineq}).

The density $Q\overline{V}(\chi(x_{\alpha}),\cdot)$ is
quasiconvex in $M^{3\times2}$ for a.e. $x\in\Omega$. Using an argument similar to that
exploited in \cite[Proposition 6]{LDR} one concludes that $Q\overline{V}%
(\chi(x_{\alpha}),\cdot)$ is quasiconvex also in $M^{3\times3}$. Thus, by the
growth condition of $Q\overline{V},$ as stated in Proposition
\ref{continuityQVbar}, the functional $v\in W^{1,p}(\Omega;\mathbb{R}%
^{3})\longmapsto\int_{\Omega}Q\overline{V}(\chi(x_{\alpha}),\nabla_{\alpha
}v(x))dx$ is sequentially weakly lower semicontinuous with respect to
$W^{1,p}$-weak topology. Hence,
\[
\liminf_{\varepsilon\rightarrow0^{+}}\int_{\Omega}Q\overline{V}(\chi
,\nabla_{\alpha}w_{\varepsilon})dx\geq2\int_{\omega}Q\overline{V}%
(\chi,\nabla_{\alpha}u)dx_{\alpha}.
\]
By the superadditivity of the lim inf we achieve the claim.

\noindent\textsl{Step two:} To prove the reverse inequality we start by observing that, fixing $\chi\in BV(\omega;\{0,1\})$, 

\noindent$ J(\chi,u)\leq\liminf_{\varepsilon\rightarrow0^{+}}%
J_{\varepsilon}(\chi,u_{\varepsilon})$ for every  $\{u_\varepsilon\}\subseteq L^p(\Omega;\mathbb R^3)$ and $ u \in W^{1,p}(\omega;\mathbb R^3)$ such that $u_{\varepsilon}\rightarrow u \hbox{ in } L^{p}(\Omega;\mathbb{R}^{3}).$

Thus it suffices to study the asymptotic behaviour with respect to the
$W^{1,p}$-weak convergence of
$$
\int_{\Omega}\left(  \chi W_{1}\left(  \nabla_{\alpha}u_{\varepsilon
}\left\vert \frac{1}{\varepsilon}\nabla_{3}u_{\varepsilon}\right.  \right)
+(1-\chi)W_{2}\left(  \nabla_{\alpha}u_{\varepsilon}\left\vert \frac
{1}{\varepsilon}\nabla_{3}u_{\varepsilon}\right.  \right)  \right)
dx-\int_{\Omega}f\cdot u_{\varepsilon}dx. 
$$
Since $\chi$ is fixed we can rewrite $\chi W_1(\cdot)+ (1-\chi)W_2(\cdot)$ as a new function with explicit dependence on $x_\alpha$.

Namely, let $W: \omega\times\mathbb{R}^{3\times3} \to \mathbb R
$ be given by 
$$
W(x_{\alpha},F):=V(\chi(x_\alpha), F)= \chi(x_{\alpha})W_{1}(F)+\left(  1-\chi(x_{\alpha})\right)
W_{2}(F), 
$$
for every $(x_{\alpha},F)\in\omega\times\mathbb{R}^{3\times3}.$

\noindent Clearly, $W$  is a Carath\'{e}odory function satisfying the growth condition
$\frac{1}{C}|F|^{p}-C\leq W(x_{\alpha},F)\leq C(1+|F|^{p}) \hbox{ for a.e. }x_{\alpha}\in\omega \hbox{ and for all } F \in\mathbb{R}^{3\times3}.$  Applying \cite[Theorem 2.3]{BaFr} to the sequence of functionals $\{G_\varepsilon\}$, where $G_\varepsilon: L^p(\Omega;\mathbb R^3 )\to [0, +\infty)$ is given by
$\displaystyle{G_\varepsilon(u):=\left\{\begin{array}{ll}\int_{\Omega}W\left(  x_{\alpha},\left( \nabla_{\alpha}u,\left\vert \frac{1}{\varepsilon}\nabla_{3}u\right.\right)\right)  dx-\int_{\Omega}f\cdot u dx &\hbox{ if }u \in W^{1,p}(\Omega;\mathbb R^3),\\
+\infty &\hbox{otherwise,}
\end{array}\right. }$

\noindent and arguing as in  \cite[Remark 3.3]{BFF}
we get 
\[
\lim_{\varepsilon \to 0^{+}}\left(  \int_{\Omega}W\left(  x_{\alpha},\left( \nabla_{\alpha}u_{\varepsilon}\left\vert \frac{1}{\varepsilon}%
\nabla_{3}u_{\varepsilon}\right.\right)\right)  dx-\int_{\Omega}f\cdot
u_{\varepsilon}dx\right)  \leq2 \int_{\omega}Q\overline{W}(x_{\alpha}%
,\nabla_{\alpha}u)dx_{\alpha}-\int_{-1}^{1}\int_{\omega}f\cdot udx_{\alpha}dx_3,\] 
where $\overline{W}:\omega\times\mathbb{R}^{3\times2}$ is defined by
$\overline{W}(x_{\alpha},\overline{F}):=\inf_{c\in\mathbb{R}^{3}}W(x_\alpha%
,(\overline{F}|c)),$
and $Q\overline{W}$ stands for the quasiconvexification of $\overline{W}$ in
the second variable.

Observing that by (\ref{Vbar})
\[
\overline{W}(x_\alpha,\overline{F})=\chi(x_\alpha)\overline{W_{1}}(\overline
{F})+(1-\chi(x_\alpha))\overline{W_{2}}(\overline{F})=\overline{V}(\chi
(x_\alpha),\overline{F}),\hbox{ and }
Q\overline{W}(x_\alpha,\overline{F})=Q\overline{V}(\chi(x_\alpha),\overline{F})
\]
for every $(x_\alpha,\overline{F})\in\omega\times\mathbb{R}^{3\times2}$, the proof is concluded.

\medskip

In the following we apply the previous analysis, with small changes, to the case where the perimeter penalization in (\ref{originalpb}) is replaced by a more general elliptic integral, such as in \cite{LK}. Namely, we consider $\Psi:\mathbb{R}%
^{3}\rightarrow\lbrack0,+\infty)$ even, continuous, positively
1-homogeneous, and satisfying (\ref{Psigrowth}).

Recall problem (\ref{originalpb2}) and observe that by $\int_{\partial E(\varepsilon)}\Psi(\nu_{E(\varepsilon
)})d\mathcal{H}^{2}$ we mean the integral with respect to the Hausdorff
measure concentrated on $\partial E(\varepsilon)$, identified with $S_{\chi
_{E(\varepsilon)}}$, with exterior normal $\nu_{E(\varepsilon)}$.

By performing the same rescaling as in (\ref{1edilation}) we obtain the
following formulation of (\ref{originalpb2}) in the fixed domain $\Omega$,%

\begin{equation}%
\begin{array}
[c]{l}%
\inf\limits_{%
\begin{array}
[c]{l}%
u\in W^{1,p}(\Omega;\mathbb{R}^{3})\\
\chi\in BV(\Omega;\{0,1\})
\end{array}
}\Big\{  {\displaystyle\int_{\Omega}}\left(  \chi W_{1}+\left(
1-\chi\right)  W_{2}\right)  \left(  \nabla_{\alpha}u\left\vert \frac
{1}{\varepsilon}\nabla_{3}u\right.  \right)  dx-\int_{\Omega}%
f\cdot u~dx \\
\qquad\qquad\qquad\qquad\quad\quad\displaystyle{+\int_{S_{\chi}}\Psi\left(
\nu_{\alpha}|\frac{1}{\varepsilon}\nu_{3}\right)  d{\cal H}^{2}}:u=0\hbox{ 
on }\partial\omega\times\left(  -1,1\right)  ,\frac{1}{\mathcal{L}^{3}%
(\Omega)}\int_{\Omega}\chi dx=\lambda\Big\},
\end{array}
\label{P2}%
\end{equation}
where $\chi$ denotes the characteristic function of $E_{\varepsilon}$, and
$\nu$ the normal to its jump set, $S_{\chi}$.

Next we will identify functions defined in $\mathbb{R}^{3}$ (or
$\mathbb{R}^{2}$) and their restrictions to $S^{2}$ (or $S^{1}$), and so the same notations will be adopted.


Let $\overline{\Psi}:\mathbb{R}^{2}\rightarrow\lbrack0,+\infty)$ be the
function given by $\overline{\Psi}(\eta):=\inf\left\{  \Psi(\eta,\xi):\xi\in\mathbb{R}\right\},$
with $\Psi:\mathbb{R}^{3}\rightarrow\lbrack0,+\infty)$ as in (\ref{Psigrowth}).

Consider the following minimum problem
\begin{equation}
\inf\limits_{%
\begin{array}
[c]{l}%
u\in W_{0}^{1,p}(\omega;\mathbb{R}^{3})\\
\chi\in BV(\omega;\{0,1\})
\end{array}
}\left\{  2\int_{\omega}Q\overline{V}(\chi,\nabla_{\alpha}u)dx_{\alpha}%
-\int_{-1}^{1}\int_{\omega}f\cdot udx_\alpha dx_3+2\int_{S_{\chi}}\overline{\Psi}%
^{\ast\ast}(\nu_{\alpha})d\mathcal{H}^{1}:\frac{1}{\mathcal{L}^{2}(\omega
)}\int_{\omega}\chi dx_\alpha=\frac{1}{2}\lambda\right\}  . \label{limitpb2}%
\end{equation}
where $\overline{\Psi}^{\ast\ast}$ denotes the convex envelope of
$\overline{\Psi}$. Namely, $\overline{\Psi}^{\ast\ast}(v):=\sup\{g:\mathbb{R}^{2} \to {\mathbb R}:g \hbox{ is convex }g(v)\leq\overline{\Psi} (v)\;\forall\,v\in\mathbb{R}^{2}\}.$

\begin{theorem}
\label{mainthm2} Consider the problems (\ref{P2}), and their
minimizers. Then the latter converge, with respect to the strong topology of
$L^{1}(\Omega;\{0,1\})\times L^{p}(\Omega;\mathbb{R}^{3})$, to the minimum of
problem (\ref{limitpb2}).

\end{theorem}

We conclude by observing that the density $\overline{\Psi}^{\ast\ast}$
satisfies all the well established properties for the lower semicontinuity of
surface integrals, such as $BV$-ellipticity (cf. \cite[Definition 5.13 and
Theorem 5.14]{AFP2}), since any continuous even function
$\phi:S^{N-1}\rightarrow\lbrack0,+\infty)$ is $BV$-elliptic if and only if its
positive 1-homogeneous extension is convex.


\section*{Acknowledgements}

\footnotesize{The authors thank Irene Fonseca for having suggested this problem and they
acknowledge the support and the kind hospitality of Departamento de
Matem\'{a}tica, Universidade de \'{E}vora and DIEII at Universit\`a di Salerno.
The work of both authors was partially supported by CIMA-UE, financed by FCT
(Funda\c{c}\~{a}o para a Ci\^{e}ncia e Tecnologia) and GNAMPA through Project ''Problemi variazionali e misure di Young nella meccanica dei Materiali Complessi''.}

\end{document}